\documentstyle{article}
\def \C{{\Bbb C}}
\def \Cbar{{\bar{\cal C}}}

\def \S{{\Sigma}}

\def \min{{\setminus}}
\def \l{{\lambda}}
\def \Z{{\Bbb Z}}

\def \ag{{\frak{g}}}

\def \O{{\cal O}}
\def \M{{\cal M}}
\def \CC{{\cal C}}
\def \XC{{{X}_{\cal C}}}
\def \MC{{{\cal M}_{\cal C}}}
\def \x{{\bf x}}
\def \k{{\kappa}}
\def \XCbar{{{X}_{\bar{\cal C}}}}
\def \MCbar{{{\cal M}_{\bar{\cal C}}}}

\def \i{{\sqrt{-1}}}

\def \g{{\gamma}}
\def \x{{\bf x}}
\def \proof{{\noindent{\it Proof.\ \ }}}
\newtheorem{th}{THEOREM}[section]

\newtheorem{prop}[th]{PROPOSITION}
\newtheorem{proposition}[th]{PROPOSITION}
\newtheorem{lem}[th]{LEMMA}

\newtheorem{cor}[th]{COROLLARY}

\title{Geometry of Moduli Spaces of Flat Bundles on Punctured Surfaces.}
\author{ Philip \ A. \ Foth$\ ^1$}
\date{}
\begin{document}
\maketitle
\input amssym.def
\begin{abstract}
For a Riemann surface with one puncture we consider moduli spaces
of flat connections such that the monodromy transformation around the puncture
belongs to a given conjugacy class with the property that a product of its
distinct eigenvalues is not equal to $1$ unless we take all of them.
We prove that these moduli spaces are smooth and their natural closures are
normal with rational singularities. \end{abstract}
%
%
%
\footnotetext[1]{The author acknowledges research support from NSF
grant DMS-9504522 for the summer of 1996.}

\section{Introduction}
\setcounter{equation}{0}

Let $X$ be a Riemann surface of the genus $g>0$ with one puncture.
We consider the moduli space of flat $GL(n, \C)$-bundles such that the
monodromy transformation around the puncture belongs to a given conjugacy class
$\CC\in SL(n, \C)$. We further assume that the class $\CC$
has property P, meaning that for the set of its
eigenvalues $(\l_1, \l_2, ..., \l_n)$
we have $\l_{i_1}\l_{i_2}\cdots \l_{i_m}\ne 1$ for any $i_1<
i_2<\dots <i_m$,\ $1\le m < n$.
We carry out all the proofs for the genus $1$ and later show that all the result
easily generalize for higher genera.

Due to the well-known correspondence between flat bundles on $X$ and
representations of its fundamental group, the problem (for $g=1$) is reduced to
the consideration of moduli spaces of pairs of matrices $(B,D)$ from $GL(n, \C)$
such that $BDB^{-1}D^{-1}\in\CC$. The main result (presented in Theorem 3.2 and
Proposition 3.4) is

\begin{th} Let $X$ be a Riemann surface with one puncture and
$\CC$ - a conjugacy class in $SL(n, \C)$ with property P.

(a) The moduli space $\MC$ of flat $\frak{gl}(n, \C)$-connections over $X$ with
the monodromy transformations around the puncture in $\CC$ is smooth.

(b) Let $\bar\CC$ is the closure of $\CC$ in $SL(n, \C)$. The variety $\MCbar$
(defined in (a) by changing $\CC$ to $\bar\CC$) is normal with rational
singularities.  \end{th}

We extend the definition of property P to the case
of orthogonal or symplectic groups to establish similar results.

We also show that if $BDB^{-1}D^{-1}$ has property P, then the pair
$(B,D)$ algebraically generates the whole group.

The moduli spaces in question are of great importance by a number of reasons.  A
theorem of Mehta-Seshadri \cite{MS} identifies two moduli spaces: the space of
unitary representations of fundamental group and the space of parabolic bundles
on $X$.  When $\CC=\exp(2\pi\i d/n)Id$ is (the class of) a central element, the
space $\MC$ is a smooth K\"{a}hler manifold and appears in algebraic geometry as
the space of holomorphic vector bundles on the closed surface of rank $n$,
degree $d$ and fixed determinant \cite{AB}.  Also those moduli spaces appear in
topological and quantum field theories; they are related to Jones-Witten
invariants (see \cite{A} for details). They are closely related to Yang-Mills
theory and geometric quantization. It is necessary to mention the results due to
Simpson \cite{S2}, which provide natural correspondence between ${\cal
D}_X$-modules, Higgs bundles and local systems on $X$ (with extra conditions).
This allows one to identify moduli spaces of those objects with the moduli
spaces we consider in the present paper.

I express deep gratitude to Jean-Luc Brylinski for helpful and valuable advices.
He kindly guided me and generously shared his knowledge.

\section{Common stabilizer of two matrices}
\setcounter{equation}{0}

We say that a matrix $C$ or a conjugacy class $\CC$ in $SL(n, \C)$ with the set
of eigenvalues $\{\l_1,\dots,\l_n\}$ has {\it property P}  if the following 
condition holds. For any $m < n$ distinct numbers $1\le i_1\le i_2\le\dots\le 
i_m\le n$ the product $\l_{i_1}\l_{i_2}\cdots \l_{i_m}$ is not equal to $1$.  
One notices that the set of matrices with property P is Zariski open in 
$SL(n,\C)$.

\begin{th} Let $B, \ D \in GL(n, \C)$ be such that $[B,D]$ satisfies property P.
Then the common stabilizer of $B$ and $D$ consists of scalar matrices only.
\end{th} \proof Let $K\in GL(n, \C)$ be non-central matrix commuting with both 
$B$ and $D$.  Let $\l$ be an eigenvalue of $K$ and let $W\subset \C^n$ be the 
kernel of $K-\l.Id$.  It follows that both $B$ and $D$ stabilize $W$. Hence the 
product of eigenvalues of $[B,D]$ which correspond to $W$ is equal to $1$.  This 
means that $W=\C^n$ and $K$ is scalar. $\bigcirc$

\

We will denote by $SL(n, \C)^2$ or $GL(n, \C)^2$ the Cartesian product of two
copies of $SL(n, \C)$ or $GL(n, \C)$ respectively.  Also denote by $$\kappa=[,]:
\ GL(n, \C)^2\to SL(n, \C)$$ the commutator map.

\

\noindent{\bf Remark.} For $n>2$ the author can prove the following statement
converse to the above Theorem. Let $\CC$ be a conjugacy class in $SL(n, \C)$.
If for any $(B,D)\in GL(n, \C)^2$ the condition $[B,D]\in\CC$ implies
$\dim Z(B,D)=1$, then $\CC$ has property P. We do not include the proof
since it is long and computational and we will not use it in the persent paper.

The situation, however, is different for $n=2$.  If we take as $\CC$ the class
of $ \pmatrix{ 1 & 1 \cr 0 & 1 } $, then for all $\{ (B,D)\in GL(2,\C)^2;\
[B,D]\in\CC\}$ the common stabilizer of $B$ and $D$ denoted by $Z(B,D)$ is the
center of $GL(n,\C)$, i.e. $\dim Z(B,D)=1$.  So, when $n=2$, if $B$ and $D$ do
not commute, then their common stabilizer is the center of $GL(n,\C)$.

\

One can think of another interpretation of property P. Let $V=\C^n$ be the
tautological representation space of $SL(n,\C)$. The spaces $\wedge^iV$ are also
naturally representation spaces for $SL(n, \C)$. Property P for the matrix $C$
means that $C$ doesn't stabilize any $\ne 0$ vector in $\wedge^iV$ for $0<i<n$.

The next proposition was first proven in \cite{Shoda}.  \begin{proposition} The 
map $\kappa$ is onto.  \end{proposition} For instance, the commutator of two 
matrices from $SL(n, \C)$ $$ \pmatrix{ 0 & e_1 & 0 & \dots & 0 \cr 0 & 0 & e_2 & 
\dots & 0 \cr . & . & . & \dots & . \cr 0 & 0 & 0 & \dots & e_{n-1} \cr e_n & 0 
& 0 & \dots & 0 } \ \ \ and \ \ \ \pmatrix{ 0 & 0 & 0 & \dots & f_n \cr f_1 & 0 
& 0 & \dots & 0 \cr 0 & f_2 & 0 & \dots & 0 \cr 0 & 0 & f_3 & \dots & 0 \cr . & 
. & . & \dots & .  } $$ can be conjugate to any semisimple element of the group.  
Also for any unipotent element $U$, $U^{-1}$ is conjugate to $U$, hence $U^2$ is 
also in the image of $\kappa$. Any unipotent element is a square of another 
unipotent element.  This proves the proposition for unipotent elements.

\section{Moduli spaces of flat bundles}
\setcounter{equation}{0}

 Let us fix a conjugacy class $\CC$ in $SL(n, \C)$.  Further we consider the 
variety of pairs and the moduli spaces of pairs of matrices with their 
commutator in $\CC$.  We denote $\XC=\{ (B,D)\in GL(n,\C)^2; [B,D]\in \CC\}$.
It is well-known that if $\CC=1$ then the variety of commuting pairs is
irreducible (\cite{MT}), but is not a smooth variety.

\begin{lem} The orbit of every element $(B,D)$, such that $[B,D] \in \CC$ is 
 closed in $\XC$ if the class $\CC$ satisfies property P.  \end{lem} \proof All 
orbits in the closure of $\XC$ have the same dimension when $\CC$ has property 
P. (Because all the conjugacy classes in the closure of $\CC$ still have 
property P, and whenever $[B,D]$ has property P the stabilizer $Z(B,D)$ 
coincides with the center of the group.) This implies that orbits are closed in 
$\XC$.  $\bigcirc$

\

Actually, we have proved a stronger result, namely that the orbit of $(B,D)$ is 
closed {\it in $GL(n, \C)^2$} if property P holds for $[B, D]$.

We consider the moduli space $$ \MC=\{ (B,D)\in GL(n,\C)^2; [B,D]\in 
 \CC\}/SL(n,\C), $$ where factoring occurs by the adjoint action of the special 
  linear group. In order to define $\MC$ properly one takes the quotient in the 
sense of the GIT (\cite{GIT}).  In the case of an affine variety $X$ there is a 
natural definition of the quotient variety $Y=X/G$. The algebra of the regular 
functions $\O(Y)$ is just the subalgebra of $\O(X)$ of $G$-invariant regular 
functions on $X$.  But when all orbits are closed, each point of the quotient 
$\MC$ correspond to an orbit in $\XC$.

Using the above lemma we see that when $\CC$ is semisimple with property P, the 
variety $\XC$ is closed affine, the quotient $\XC/SL(n,\C)$ exists and its 
points correspond exactly to the orbits. We can construct a nice quotient when 
$\CC$ is not semisimple, but still satisfies property P.  Let $\Cbar$ be the 
closure of $\CC$ in $SL(n,\C)$ and $\XCbar=\{(B,D)\in 
GL(n,\C)^2;[B,D]\in\Cbar\}$. It is affine so the quotient $\MCbar$ by $SL(n,\C)$ 
exists.  Since all $SL(n,\C)$-orbits have the same dimension, points of $\MCbar$ 
again correspond exactly to orbits. It is well-known that $\Cbar\min\CC$ is the 
union of {\it finite} number of conjugacy classes $\CC_i$. So we have the 
corresponding finite set of closed subvarieties $\M_i$ in $\MCbar$.  We define 
now the algebraic variety $\MC$ which is the complement of the union 
$\cup_i\M_i$ in $\MCbar$.

\

{\it Example.} Here we consider the simple but important case $n=2$.  We have 
five types of conjugacy classes in $SL(2, \C)$ (here we mention just a 
representative of each):  $I=\pmatrix{ 1 & 0 \cr 0 & 1}$ ($\dim (I)=0$),
$-I=\pmatrix{ -1 & 0 \cr 0 & -1}$ ($\dim (-I)=0$), $R_2=$ the class of 
$\pmatrix{ 1 & 1 \cr 0 & 1}$, ($\dim (R_2)=2$), $R_e=$ the class of $\pmatrix{ 
-1 & 1 \cr 0 & -1}$, ($\dim (R_e)=2$), $R_{\l}=$  the class of $\pmatrix{ \l & 0 
\cr 0 & \l^{-1}}$, $\l^2\ne 1$, ($\dim (R_{\l})=2$).  One has $$ GL(2, 
\C)^2=X_I\cup X_{-I}\cup X_{R_2}\cup X_{R_e}\cup \bigcup_{\l^2\ne 1}X_{R_{\l}}.  
$$ Each $\XC$ except for $X_I$ is a connected smooth variety. Their dimensions 
are $6,\ 5,\ 7,\ 7,\ 7$ respectively.  The varieties $\XC$ corresponding to
semisimple classes are closed.  Of course, $X_I$ and $X_{-I}$ lie in the closure 
of $X_{R_2}$ and $X_{R_e}$ respectively. Also, they both are limit varieties of
$X_{R_{\l}}$ for $\l\to\pm 1$, so they both lie in the closure of $\cup_{\l^2\ne 
1}X_{R_{\l}}$. As it was mentioned above, $X_I$ is an irreducible variety. Also 
it has singularities only in codimension 2.

The space ${\cal M}_{\bar{R_e}}$ identifies with the space of pairs of matrices 
$(B,D)$ from $GL(2,\C)$ of the form $B=\pmatrix{-\l & x \cr 0 & \l}$,
$D=\pmatrix{y & a \cr a & 0}$ and it has dimension $4$. It is isomorphic to 
$\C^*\times \C^*\times \C^2$ as an algebraic variety. Clearly ${\cal 
M}_{\bar{R_e}}={\cal M}_{R_e}\cup {\cal M}_{-I}$. The space ${\cal M}_{-I}$ 
identifies with the space of pairs of matrices $(B,D)$ from $GL(2,\C)$ of the 
form $B=\pmatrix{-\l & 0 \cr 0 & \l}$, $D=\pmatrix{0 & a \cr a & 0}$, which is 
isomorphic to $\C^*\times \C^*$ and has dimension $2$. (So, $\MC$ is isomorphic 
to to $\C^*\times \C^*\times \C^2\min {0}$, and hence is not affine.)

In general, property P for the class $\CC$ implies that the affine variety 
$\MCbar$ is stratified by the smooth locally closed subvarieties $\MC$ and 
${\cal M}_{\CC_i}$ and the codimension of ${\cal M}_{\CC_i}$ in $\MCbar$ is 
equal to the codimension of $\CC_i$ in $\Cbar$. This means, first of all, that 
all the singularities of $\MCbar$ are in codimension at least $2$.  \begin{th} 
If $\CC$ satisfies property P then $\MCbar$ is a normal affine algebraic 
variety.  It is Cohen-Macaulay and has rational singularities.  \end{th} Before 
we prove this assertion, we exhibit an auxiliary result.  \begin{lem} If 
$[B,D]=A$ and $\dim Z(B,D)=1$ then the differential map $d\k: 
T_{(B,D)}GL(n,\C)^2\to T_ASL(n,\C)$ is surjective${}^1$.  \end{lem}
\footnotetext[1]{Jean-Luc Brylinski noticed that it should be true in every
characteristic.} \proof We identify the first tangent space with 
$\frak{gl}(n,\C)^2$ and the second one with $\frak{sl}(n,\C)$ via the 
corresponding left multiplications.  Computations show that with these 
identifications $d\k$ sends $(x,y)\in \frak{gl}(n,\C)^2$ to
$DB(D^{-1}xD-x+y-B^{-1}yB)B^{-1}D^{-1}$. So, it is enough to show that
$R(AdD-1)+R(AdB-1)=\frak{sl}(n,\C)$, where $R(L)$ denotes the range of a linear 
operator $L$. With respect to the bilinear form $Tr(XY)$ one has
$R(AdB-1)=Ker(AdB-1)^{\perp}$. We notice that $Ker(AdB-1)=z(B)$ - the 
centralizer of $B$ in the Lie algebra.  Now the condition of the lemma implies 
that $$ R(AdD-1)+R(AdB-1)=Ker(AdB-1)^{\perp} + Ker(AdD-1)^{\perp}= $$ $$ 
=(z(B)\cap z(D))^{\perp}=z(B,D)^{\perp}=\frak{sl}(n,\C). \  \  \ \bigcirc $$

\noindent{\it Proof of the theorem.} Let us consider the morphism of algebraic 
varieties $\k : \XCbar\to\Cbar$. We notice that every conjugacy class in the 
closure of a conjugacy class with property P also satisfies property P.  To see 
that $\k$ is actually a smooth morphism, one notices that smoothness is 
preserved by base extensions. Let $U$ be the open set of matrices with property 
P and $U^{(2)}\subset GL(n, \C)^2$ its preimage under $\k$. We saw in Lemma 3.3 
that $\k : U^{(2)}\to U$ is smooth. For a class $\CC$ with property P we have 
the inclusion $\Cbar\hookrightarrow U$.  Now we make the base change : 
$U^{(2)}\times_U \Cbar\to\Cbar$.

We conclude that $\k : \XCbar\to\Cbar$ is a smooth morphism.  Now we invoke the 
theorem of Kraft and Procesi (\cite{KP}) which tells us that $\Cbar$ is normal, 
Cohen-Macaulay with rational singularities.  As a consequence we obtain the fact 
that $\XCbar$ has rational singularities. Now we use the theorem of Boutot 
\cite{Bou}, which implies that $\MCbar$ has also rational singularities. In 
particular, $\MCbar$ is normal and Cohen-Macaulay.  $\bigcirc$

\

Actually, we have proved that every connected component of $\XCbar$
(and $\MCbar$) is normal. But it seems likely that those varieties
are connected. As far as we know this is an open problem.

It is very possible that whenever a class $\CC$ in the image of $\k$
is not the identity, we can find a pair $(B,D)\in\XC$ with
trivial $Z(B,D)$. If it would be so, one could reprove the results of Shoda
and Ree (Propositions 2.2 and 5.2) in a nice algebraic way as follows.
Let $X_I\subset G\times G$ be the variety of commuting matrices. The above
result on the surjectivity of $d\k$ implies that the image of $\k$
restricted to $(G\times G)\min X_I$ is open in $G\min \{I\}$. We know
from previous explicit constructions that every semisimple or unipotent
element is in the image of $\k$. One needs only to remark
that each conjugacy class in $G\min \{ I\}$ is either unipotent
or contains in its closure a semisimple element $\ne I$.

\begin{proposition} If $\CC$ has property P then $\MC$ is a smooth algebraic 
variety.  \end{proposition} \proof The statement follows from the above theorem, 
because the fibers of the map $\k: \XCbar\to\Cbar$ are smooth varieties. 
$\bigcirc$

\

By a procedure similar to the one we described above, it is possible to define
the space $\MC$ for {\it any} conjugacy class $\CC$. (First, we define $\MCbar$
and then using its natural stratification we throw out irrelevant pieces.) The
following lemma which was pointed out to me by J.-L. Brylinski calculates the
dimension of the tangent space to $\MC$ at the class $(B,D)$.  A point
$(B,D)\in\XC$ is called a general point if $\XC$ is smooth at $(B,D)$ and all
stabilizers of elements in some neighbourhood of $(B,D)$ in $\XC$ form a smooth
group bundle.

\begin{lem} Let $\CC$ be a conjugacy class in $SL(n, \C)$. For a
general point $(B, D)\in\XC$ one has $\dim (T_{(B,D)}\MC)=\dim \CC +2\dim
Z(B,D).$ \end{lem} \proof The loops $a$ and $b$ generate freely the fundamental
group of the elliptic curve with one puncture $\S\setminus \{O\}$.  The
corresponding monodromy transformations $B$ and $D$ define the local system $V$
on ${\S}$ of the dimension $n$.  Let us define a manifold $X$ as consisting of
all conjugacy classes of homomorphisms $\rho: \pi_1(\S\min\{O\})\to SL(n,\C)$
such that the image of $\rho$ is Zariski dense in $SL(n,\C)$. (So that $\MC$ is
a subspace of $X$.) The tangent space to $X$ in the class of $\rho$ identifies,
by a well-known theorem of A. Weil, with the group cohomology
$H^1(\pi_1(\S\min\{O\}),\ag)$, where $\ag=\frak{sl}(n, \C)$ is a
$\pi_1(\S\min\{O\})$-module via the adjoint action followed by $\rho$.  For any
$A\in SL(n,\C)$ we identify the tangent space $T_ASL(n, \C)$ to $\ag$ via the
action of the left translation by $A$.  The tangent space to the conjugacy class
$\CC_A$ of $A$ is the subspace of $\ag$ given as the range of $Ad(A)-1:
\ag\to\ag$. If $\Gamma$ is the cyclic subgroup generated by $\g=aba^{-1}b^{-1}$
then the cohomology group $H^1(\Gamma,\ag)$ is the cokernel of the map
$Ad(A)-1$. So the tangent space $T_{[\rho ]}\MC$ identifies with the subspace
$Ker(H^1(\pi_1(\S\min\{O\}),\ag)\to H^1(\Gamma,\ag))$ of $T_{[\rho
]}X=H^1(\pi_1(\S\min\{O\}),\ag)$. Thus, \begin{equation}
T_{(B,D)}\MC=Ker[H^1({\S}\setminus \{O\}, End(V))\to H^1(D^*, End(V))],
\end{equation} where $D^*$ is a small disk around the puncture. (It is the same
as $$ Ker[H^1(\pi_1({\S}\setminus \{O\}), End(V))\to H^1(\Z , End(V))], $$
because $\pi_1(D^*)=\Z$.) The Euler characteristic of the punctured elliptic
curve is $\chi({\S}\setminus \{O\})=-1$ and whence, $$ \dim [H^0({\S}\setminus
\{O\}, End(V))]-\dim [H^1({\S}\setminus \{O\}, End(V))]= $$ $$
=n^2\chi({\S}\setminus \{O\})=-n^2.  $$ But $H^0({\S}\setminus \{O\},
End(V))=End(V)^{B,D}=Z(B,D)$ in $GL(n, \C)$.  So, \break $\dim H^1({\S}\setminus
\{O\}, End(V))=n^2+\dim   Z(B, D)$.  Also we notice that the map in the equation
(3.1) is onto due to the exact sequence${}^1$ \footnotetext[1]{The fact that
$T_{(B,D)}\MC = Im(H^1_c(\S\setminus\{0\}, End(V))\to H^1(\S\setminus\{0\},
End(V)))$ was also used in \cite{BG}.} $$ H^1({\S}\min \{O\}, End(V))\to
H^1(D^*, End(V))\to H^2_c({\S}\min\{O\}, End(V)) \to $$ $$ \to
H^2({\S}\min\{0\}, End(V))\to 0.  $$ The group $H^2 ({\S}\min \{ 0\}, End(V))$
is $0$ for dimension reasons.  Besides, one sees that the group $H^2_c ({\S}
\min \{ 0\}, End (V))$ is dual to $H^0(({\S}\min \{ 0\}, End(V))$, because the
local system $End(V)$ is self-dual. The group $H^0(({\S}\min \{ 0\}, End(V))$ is
the group of matrices which commute with $B$ and $D$. It follows that the image
of the linear map $H^1({\S}\min \{ 0\}, End(V)) \to H^1(D^*, End(V))$ has
dimension equal to $$\dim  H^1(D^*,End(V)) - \dim Z(B,D) = \dim Z([B,D]) - \dim
Z(B,D).$$ Therefore the dimension of the kernel of this same map is equal to
$$\dim H^1({\S}\min \{ 0\},End(V))-\dim Z([B,D])+\dim Z(B,D)$$ which is equal to
$$n^2+\dim Z(B,D)-\dim Z([B,D])+\dim Z(B,D)=\dim (\CC)+2\dim Z(B,D).$$ Here we
used the fact that $n^2 = \dim GL(n, \C) = \dim (\CC) + \dim Z(g)$, $g\in \CC$.
$\bigcirc$

\

First, we notice that the space $\MC$ has even dimension, (the smooth loci are
actually known to be hyper-K\"{a}hler manifolds). Also we have \begin{cor} $\dim
(\XC)=n^2+\dim \CC +\dim Z(B,D)$, where $(B,D)$ is a generic element of $\XC$.
\end{cor}

\section{On subalgebras generated by pairs}
\setcounter{equation}{0}

Let $p$ be an integer, $p>1$, and let $G^p=\underbrace{G\times\dots\times G}_p$,
where $G$ is a reductive algebraic group (over $\C$) with Lie algebra $\ag$. Let
$\x=(x_1,\dots,x_p)\in G^p$, $G.\x$ - its orbit, and $A(\x)$ - the algebraic
subgroup of $G$ generated by the set $\{x_1,\dots,x_p\}$. (So $A(\x)$ is the
Zariski closure of the subgroup of $G$ in the abstract sense generated by the 
set $\{x_1,\dots,x_p\}$.) Let also $\pi: G^p\to G^p/G$ stand for the quotient 
morphism. Following Richardson \cite{Rich} we call $\x$ a {\it semisimple 
$p$-tuple} if $A(\x)$ is linearly reductive. (Since we are in characteristic 
zero this it is equivalent to reductive.) We cite from \cite{Rich} the following
\begin{th} The orbit $G.\x$ is closed if and only if $\x$ is semisimple.
\end{th} This allows us to apply our knowledge to the situation of the group
$GL(n,\C)$ and $p=2$.  \begin{lem} If $\CC$ satisfies property P and 
$[B,D]\in\CC$ then the algebraic subgroup $A(B,D)$ of $GL(n,\C)$ is reductive.
\end{lem} \proof We saw before (sections 2 and 3) that property P for $\CC$
implies that the orbit of the element $(B,D)$ is closed. Now the above theorem 
finishes the job. $\bigcirc$

\begin{prop} If $\CC$ satisfies property P and $[B,D]\in\CC$, then
$A(B,D)=GL(n,\C)$.  \end{prop} \proof Let $\frak a\subset \frak{gl}(n,\C)$ be 
the Lie algebra of $A(B,D)$.  The lemma above implies that $\frak a$ is 
reductive and its centralizer inside $\frak{gl}(n,\C)$ is the center of the 
algebra. Now we apply the double commutant theorem to obtain that $\frak a$ 
actually coincides with $\frak{gl}(n,\C)$.  $\bigcirc$

\section{Symplectic and orthogonal groups} 
\setcounter{equation}{0} 

Here we will formulate property P for orthogonal and symplectic groups.  It 
turns out that many of the results that we proved remain valid for other 
semisimple algebraic groups.  Let $G$ be either $SO(2n, \C)$, $SO(2n+1, \C)$, or 
$Sp(2n, \C)$. The group $G$ naturally acts on $\C^{2n(+1)}$ preserving the 
bilinear form $\langle , \rangle$.  It is well-known that if $\l$ is an 
eigenvalue of $A\in G$ then $\l^{-1}$ is an eigenvalue too of the same 
multiplicity and partition.

Let $\l_1^{\pm 1}, \l_2^{\pm 1}, ..., \l_n^{\pm 1}$ (and $1$ in the case of 
$SO(2n+1, \C)$) be the set of eigenvalues of a conjugacy class $\CC\subset G$.
We say that $\CC$ has {\it property P} if no product of its eigenvalues of the 
form $\prod_{j\in S}\l_j^{e_j}$ is equal to one, where $S$ is non-empty subset 
of $\{ 1,2,...,n\}$ and $e_j=\pm 1$.  For an element $A\in \CC$ this is 
equivalent to the condition that for any isotropic subspace $V$ preserved by $A$ 
the product of eigenvalues of $A$ in $V$ is not equal to $1$.

\

\noindent{\it Remark.${}^1$} One can formulate property P for an element
$C$ for any group $G$ as follows. Let $C=C_sC_u$ be a Jordan decomposition into
a product of commuting unipotent and semisimple elements and let $\C^r$ be the
standard representation space of $G$. Consider $V=\oplus_{i=0}^r\wedge^i\C^r$,
which is naturally a representation space of $G$ too. Now we say that $C$ (or
its conjugacy class in $G$) has property P if the following two stable subspaces
of $V$ have the same dimension: $\dim(V^{C_s})=\dim(V^T)$, where
$T$ is a maximal torus of $G$.

\
\footnotetext[1]{The idea of this nice remark is due to Ranee Brylinski.}

Suppose that an element $K\in G$ belongs to the common stabilizer of elements 
$B$ and $D$. Also we make an assumption that $K^2\ne I$ - the identity matrix. 
We will show that we may always find an isotropic subspace $V$ preserved by both 
$B$ and $D$. It means in turn that the product of eigenvalues of $[B,D]$ in $V$ 
is equal to $1$.  At first, we consider the case when $K$ has an eigenvalue 
$\l\ne\pm 1$.  Here we may take as $V$ the kernel of $K-\l.Id$.

It remains to assume that $K$ has only $\pm 1$ as the set of eigenvalues.  The 
condition that $K^2\ne I$ implies that there exists a vector $y\in \C^{2n(+1)}$ 
such that $y$ belongs to the generalized eigenspace of $1$ or $-1$ and $Ky\ne\pm 
y$.  It is equivalent to the condition that there is an eigenvalue $\mu$ of $K$ 
(which is, of course, $1$ or $-1$), and a $K$-irreducible and $K$-invariant 
subspace $Y$ such that $\dim Y >1$.  Let us consider the subspace 
$W=Ker(K-\mu.Id)$. Define the subspace $V\subset W$ consisting of all $v\in W$ 
such that the equation $K^{-1}x=\mu x+v$ has a solution.  One notices that $V$ 
is the intersection of $Ker(K-\mu)$ with $Im(K^{-1}-\mu)$. Then for $v\in V$ we 
have $$ \langle v,v\rangle =\langle v, K^{-1}x-\mu x\rangle =\langle Kv-\mu 
v,x\rangle =0.  $$ Moreover, if $B$ commutes with $K$, then $V$ is also 
$B$-invariant.  The fact that $V\ne 0$ follows from the choice of $\mu$.

Now we will simply rewrite some of the statements which we proved before for the
group $SL(n, \C)$ and which continue to be true in the case when $G=SO(2n, \C)$, 
$SO(2n+1, \C)$, or $Sp(2n, \C)$. We just saw that \begin{th} Let $B, D\in G$ be 
such that $[B,D]$ satisfies property P.  Then the common stabilizer of $B$ and 
$D$ is finite of exponent $2$.  \end{th} Let $G^2$ stand for the cartesian 
product of two copies of $G$ and let $\k=[,]: G^2\to G$ be the commutator map.  
The next proposition is quoted from \cite{Ree}:  \begin{prop} The map $\k$ is 
onto.  \end{prop} Analogously to Lemma 2.3 one may prove the following
\begin{lem} For $n>2$ and a conjugacy class $\CC\subset G$ the following two 
conditions are equivalent:

(i)\ \ \ $\CC$ has property P.

(ii)\ For any $(B,D)\in G^2$ such that $[B,D]\in \CC$ we have $\dim Z(B,D)=0$.
\end{lem}

Also we have \begin{lem} If $\CC$ has property P and $[B,D]\in \CC$ then the 
orbit of the element $(B,D)\in G^2$ is closed.  \end{lem}

Just as in section 3, we may define the space $\XC$ and the moduli space $\MC$.
\begin{prop} (a) $\dim(\MC)=\dim\CC+2\dim Z(B,D)$  and $\dim 
\XC=\dim(G)+\dim\CC+\dim Z(B,D)$, for a generic element $(B,D)\in \XC$.

(b) If $\CC$ has property P, then $\MC$ has at worst quotient singularities.
\end{prop}

Unfortunately, we cannot say anything about the normality of $\MCbar$,
because the closures of conjugacy classes in $G$ are not always normal
(cf. \cite{KP}).

The result from the section 4 about the algebraic subgroup $A(B,D)$ generated by
a pair $(B,D)\in G^2$ still generalizes:  \begin{prop} If the commutator $[B,D]$
belongs to a conjugacy class $\CC$ with property P, then $A(B,D)=G$.  \end{prop}

\noindent{\bf Remark.}  We saw that the questions of normality of variety
$\MCbar$ and rationality of its singularities depend only on the geometry of
$\Cbar$ if the class $\CC$ has propetry P. In the case of orthogonal or
symplectic group there occur different phenomena which are described by Kraft
and Procesi in \cite{KrP}.  Using their results one can easily prove the analog
of the part {\it (b)} of the above theorem for the group $SO(n, \C)$ or $Sp(2n,
\C)$ according to whenever $\Cbar$ is normal.

\section{Generalizations for $p$-tuples.}
\setcounter{equation}{0}

Here we briefly discuss how results of the paper can be extended to the case of 
$p$-tuples.

Let us consider the operation $\k:=[,]: G^p\to G$, $$ [A_1, A_2, ..., 
A_p]=A_1A_2\cdots A_pA_1^{-1}A_2^{-1}\cdots A_p^{-1}.  $$ It is clear what 
property P means in this case. It is still true that if \break $[A_1, A_2, ..., 
A_p]\in\CC$ has property P, then the stabilizer $Z(\x)$ (where $\x :=\break(A_1, 
A_2, ... A_p)$) is the center of $G$ in the case $G=SL(n, \C)$ and is finite of 
exponent $2$ when $G$ is orthogonal or symplectic group.  Taking $A_3=A_4=\cdots 
=A_p=1$ we see that this fact actually characterizes the class $\CC$ and that 
the map $\k$ is onto. If $[A_1, A_2, ..., A_p]$ has property P, then its orbit 
is closed.  As before, we can define the spaces $\XC$ and $\MC$. They are of 
dimensions $(p-1)\dim(G)+\dim\CC+\dim Z(\x)$ and $(p-2)\dim(G)+\dim\CC+2\dim 
Z(\x)$. When $G=GL(n, \C)$ and the class $\CC$ has property P, both $\XC$ and 
$\MC$ are smooth varieties.  Notice, that in the application to Riemann surfaces 
$p$ is even and so is the dimension of $\MC$.

In the case $G=SL(n, \C)$ the same method as in section 3 proves the normality 
of $\MCbar$. We also deduce that if $\x$ has property P, then the set $\{ A_1, 
A_2, ..., A_p\}$ generates the whole group $G$.

\

\noindent{\bf Remark.} Let ${\bar X}$ be a compact oriented surface of genus
$g>0$, let $s_1, ..., s_k\in {\bar X}$ and let $X={\bar X}\min\{ s_1, ...,
s_k\}$ be the punctured surface. We consider local systems of rank $n$ on $X$.
Assume that $a_1, a_2, ..., a_{2p}$ are loops in $X$ generating the fundamental 
group $\pi_1(\bar X, x_0)$.  Also let $\g_i$ be a loop which comes from $x_0$, 
goes once around $s_i$ counterclockwise and then goes back to $x_0$. One can 
take them in such a way that $\g_1\g_2\cdots \g_k = [a_1, a_2, ..., a_{2p}]$.  
Let us have $k$ matrices $C_1, ...,C_k$ defined up to simultaneous conjugation, 
such that $C_i$ is the monodromy transformations of a local system $V$ on $X$ 
corresponding to $\g_i$.  So we have the corresponding matrix equation:
$C_1C_2\cdots C_k=[A_1, A_2, ..., A_{2p}]$.  Here $A_1, ..., A_{2p}$ are the 
monodromy transformations corresponding to the generators of the fundamental 
group of the compactified curve. (The local system $V$ is comletely determined
by the $2p+k-1$ matrices $C_1, C_2, ..., C_{k-1},A_1, A_2, ..., A_{2p}$.) Under 
the assumption that $\det(C_1\cdots C_k)=1$ Proposition 2.2 (and 5.2) shows that 
we always can find a solution of this equation. When the genus is zero this is 
not always the case (see \cite{Simmat} for details).

 \end{document}